\documentstyle[a4,12pt]{article}
\begin{document}
\begin{center}
{\bf Differential geometry of surfaces and  Heisenberg ferromagnets }
\end{center}
\begin{center}
H.S.Serikbaev, Kur.Myrzakul,  F.K.Rahimov
\end{center}
\begin{center}
{\it Institute of Physics and Technology,
480082, Alma-Ata-82, Kazakhstan}
\end{center}
\begin{abstract}
The relation between differential geometry of  surfaces and some Heisenberg
ferromagnet models is  considered. 
\end{abstract}
%\tableofcontents
\section{Introduction}
The deep interrelation between many partial differential equations 
(PDE) of
the classical differential geometry (DG) of surfaces
and  modern soliton equations is well established now [2-8]. Generally speaking,
the interaction between DG and PDE has been studied
since the 19$^{th}$ century. For example, it can be found in the classical
works of Lie, Darboux, Goursat, Bianchi, Backlund, E.Cartan. In particular,
there are many very interesting PDE originating from DG.
The famous sine-Gordon equation, for example, first appeared in DG.  Some
of the other soliton equations have been known in DG for a long time [9-11].

One of interesting subclass of soliton equations, useful both from the mathematical
and physical point of view, is Heisenberg ferromagnets (HF) or 
spin systems (SS) [1,12-27].  
In the present paper, we will construct some classes of  surfaces corresponding
several  SS in two dimensions.

\section{Some fundamental formulas on the theory of surfaces}
Let $M^{2}$ be a surface in the 3-dimensional space $E^{3}$ parametrized
by the coordinates $x, y$ and defined by the position vector ${\bf r}(x,y)$.
The unit normal vector to the surface
 defines as usual
$$
{\bf n}=\frac{{\bf r}_x\wedge {\bf r}_y}{|{\bf r}_x\wedge {\bf r}_y|}.\eqno(1)
$$
Now we cite a few formulas for the  tangent vectors which induce some surfaces.
\subsection{The Rodrigues formula}
The simplest  example is  the  Rodrigues formula (RF) [11]
$$
{\bf r}_x=-\rho _1{\bf N}_x, \quad {\bf r}_y=-\rho _2{\bf N}_y \eqno (2)
$$
where  ${\bf N}$ is in general some vector-function, $\rho_{i}$ is the scalar function.
\subsection{The Lelieuvre formula}
Our second example is  the Lelieuvre formula (LF) [4, 7, 9-11, 19-23]
$$
{\bf r}_x=\rho {\bf N}_x\wedge {\bf N}, \quad {\bf r}_y=\rho{\bf N}\wedge{\bf N}_y. 
\eqno (3)
$$

\subsection{The Schief formula}
One of interesting generalizations of the LF (3) is  the Schief formula (SF)
[11]
$$
{\bf r}_x=\rho {\bf N}_x\wedge {\bf N}+\mu {\bf N}_x, \quad {\bf r}_y=
\rho{\bf N}\wedge{\bf N}_y+\mu {\bf N}_y. 
\eqno (4)
$$
\subsection{The Myrzakulov formula}

At last, we  consider   the Myrzakulov formula (MF) [27-28]
$$
{\bf r}_{x}=a_{1}{\bf N}\wedge{\bf N}_{x}+a_{2}{\bf N}\wedge{\bf N}_{y}+
a_{3}{\bf N}_{x}+a_{4}{\bf N}_{y}+a_{5}{\bf N}
\eqno(5a)
$$
$$
{\bf r}_{y}=b_{1}{\bf N}\wedge{\bf N}_{x}+b_{2}{\bf N}\wedge{\bf N}_{y}+
b_{3}{\bf N}_{x}+b_{4}{\bf N}_{y}+b_{5}{\bf N}
\eqno(5b)
$$
where  ${\bf N}$ is  some vector function, $a_{i},
 b_{i}$ are in general some real functions. We note that the RF (2), the LF (3) and
 the SF (4) are the particular cases of the MF (5).
\section{Surfaces and ${\bf N}$-systems}
 The compatibility condition of the
equations (5) gives [28]
$$
(a_{1y}-b_{1x}){\bf N}\wedge{\bf N}_{x}+(a_{2y}-b_{2x}){\bf N}\wedge{\bf N}_{y}+
(a_{1}+b_{2}){\bf N}_{y}\wedge{\bf N}_{x}+
$$
$$
+{\bf N}\wedge[a_{2}{\bf N}_{yy}+
(a_{1}-b_{2}){\bf N}_{xy}-b_{1}{\bf N}_{xx}]+(a_{5y}-b_{5x}){\bf N}+
(a_{3y}-b_{5}-b_{3x}){\bf N}_{x}+
$$
$$
+(a_{5}+a_{4y}-b_{4x}){\bf N}_{y}+[(a_{3}-b_{4}){\bf N}_{xy}+a_{4}{\bf N}_{yy}
-b_{3}{\bf N}_{xx}]=0.
\eqno(6a)
$$
Hence we obtain
$$
b_{5x}-a_{5y}=\frac{1}{{\bf N}\cdot {\bf N}}\{(a_{3y}-b_{5}-
b_{3x}){\bf N}\cdot{\bf N}_{x}+
+(a_{5}+a_{4y}-b_{4x}){\bf N}\cdot{\bf N}_{y}
$$
$$
+(a_{1}+b_{2}){\bf N}\cdot ({\bf N}_{y}\wedge{\bf N}_{x})+
{\bf N}\cdot [(a_{3}-b_{4}){\bf N}_{xy}+a_{4}{\bf N}_{yy}
-b_{3}{\bf N}_{xx}]\}.
\eqno(6b)
$$

Equations (6) we call the ${\bf N}$-equations or ${\bf N}$-systems
[28].
\section{Surfaces and ${\bf n}$-systems or spin systems}
Let us now we consider the case when
$$
{\bf N}={\bf n}\equiv {\bf S} \eqno(7)
$$
where ${\bf S}$ is the spin vector
$
{\bf S}=(S_{1}, S_{2}, S_{3}), \quad 
{\bf S}^{2}=S^{2}_{1}+ S^{2}_{2}+ S^{2}_{3}=1.$
In this case, the equations (6)  we call the ${\bf n}$-equations or spin systems
(SS).
Let us consider  the MF [28]
$$
{\bf r}_{x}=a_{1}{\bf S}\wedge{\bf S}_{x}+a_{2}{\bf S}\wedge{\bf S}_{y}+
a_{3}{\bf S}_{x}+a_{4}{\bf S}_{y}+a_{5}{\bf S}
\eqno(8a)
$$
$$
{\bf r}_{y}=b_{1}{\bf S}\wedge{\bf S}_{x}+b_{2}{\bf S}\wedge{\bf S}_{y}+
b_{3}{\bf S}_{x}+b_{4}{\bf S}_{y}+b_{5}{\bf S}
\eqno(8b)
$$
which is the spin form of the formulas (5). At the same time, after (7) the equations 
(6) take the form [28]
$$
(a_{1y}-b_{1x}){\bf S}\wedge{\bf S}_{x}+(a_{2y}-b_{2x}){\bf S}\wedge{\bf S}_{y}+
(a_{1}+b_{2}){\bf S}_{y}\wedge{\bf S}_{x}+
$$
$$
+{\bf S}\wedge[a_{2}{\bf S}_{yy}+
(a_{1}-b_{2}){\bf S}_{xy}-b_{1}{\bf S}_{xx}]+(a_{5y}-b_{5x}){\bf S}+
(a_{3y}-b_{5}-b_{3x}){\bf S}_{x}+
$$
$$
+(a_{5}+a_{4y}-b_{4x}){\bf S}_{y}+[(a_{3}-b_{4}){\bf S}_{xy}+a_{4}{\bf S}_{yy}
-b_{3}{\bf S}_{xx}]=0
\eqno(9a)
$$
$$
b_{5x}-a_{5y}=(a_{1}+b_{2}){\bf S}\cdot ({\bf S}_{y}\wedge{\bf S}_{x})+
{\bf S}\cdot [(a_{3}-b_{4}){\bf S}_{xy}+a_{4}{\bf S}_{yy}
-b_{3}{\bf S}_{xx}].
\eqno(9b)
$$
From (9a) as the particular case as 
$$
a_{1y}-b_{1x}=a_{2y}-b_{2x}=b_{4}-a_{3}=a_{4}=b_{3}=0
$$
we get the equation
$$
(a_{1}+b_{2}){\bf S}_{y}\wedge{\bf S}_{x}+{\bf S}\wedge[a_{2}{\bf S}_{yy}+
(a_{1}-b_{2}){\bf S}_{xy}-
b_{1}{\bf S}_{xx}]+(a_{5y}-b_{5x}){\bf S}+
$$
$$
(a_{3y}-b_{5}){\bf S}_{x}+
(a_{5}-a_{3x}){\bf S}_{y}=0.
\eqno(10)
$$
One of interesting reduction of this equation is  the stationary Myrzakulov
XIII (M-XIII) equation [28] 
$$
{\bf S}\wedge[a_{2}{\bf S}_{yy}+
(a_{1}-b_{2}){\bf S}_{xy}-b_{1}{\bf S}_{xx}]+
(a_{3y}-b_{5}){\bf S}_{x}+
(a_{5}-a_{3x}){\bf S}_{y}=0
\eqno(11a)
$$
$$
a_{5y}-b_{5x}=(a_{1}+b_{2}){\bf S}\cdot ({\bf S}_{x}\wedge{\bf S}_{y}).
\eqno(11b)
$$
It is the stationary case of the following  M-XIII equation
$$
{\bf S}_{t}={\bf S}\wedge[a_{2}{\bf S}_{yy}+
(a_{1}-b_{2}){\bf S}_{xy}-b_{1}{\bf S}_{xx}]+
(a_{3y}-b_{5}){\bf S}_{x}+
(a_{5}-a_{3x}){\bf S}_{y}
\eqno(12a)
$$
$$
a_{5y}-b_{5x}=(a_{1}+b_{2}){\bf S}\cdot ({\bf S}_{x}\wedge{\bf S}_{y}).
\eqno(12b)
$$
\subsection{Heisenberg ferromagnet}
Let us consider the following particular case of the system (8)
$$
a_{1}=a_{2}=a_{3}=a_{4}=b_{2}=b_{3}= b_{4}=b_{5}=0, 
\quad a_{5}=b_{1}=1. \eqno(13)
$$
So we have (see, e.g., [5-6])
$$
{\bf r}_{x}={\bf S}, \quad 
{\bf r}_{y}={\bf S}\wedge{\bf S}_{x}.
\eqno(14)
$$

Then the system (10) takes the form
$$
{\bf S}_{y}={\bf S}\wedge{\bf S}_{xx}
\eqno(15)
$$
that is the HF equation.
\subsection{Stationary Landau-Lifshitz equation in two dimensions}
Let us consider the following particular case of the system (8)
$$
{\bf r}_{x}={\bf S}\wedge{\bf S}_{y}, \quad
{\bf r}_{y}=-{\bf S}\wedge{\bf S}_{x}.
\eqno(16)
$$
Then the system (9) takes the form (see, e.g., the ref. [11])
$$
{\bf S}\wedge({\bf S}_{yy}+
{\bf S}_{xx})=0
\eqno(17)
$$
that is the  stationary Landau-Lifshitz equation (LLE) which follows from
the (2+1)-dimensional LLE
$$
{\bf S}_{t}={\bf S}\wedge({\bf S}_{yy}+
{\bf S}_{xx}).\eqno(18)
$$

\subsection{Stationary M-XIIIA equation}
Now we take the following  case of the system (8)
$$
a_{5}=a_{3x}+\phi_{x}, \quad b_{5}=a_{3y}-\phi_{y}. \eqno(19)
$$
In this case  we have the MF [28]
$$
{\bf r}_{x}=a_{1}{\bf S}\wedge{\bf S}_{x}+a_{2}{\bf S}\wedge{\bf S}_{y}+
a_{3}{\bf S}_{x}+(a_{3x}+\phi_{x}){\bf S}
\eqno(20a)
$$
$$
{\bf r}_{y}=b_{1}{\bf S}\wedge{\bf S}_{x}+b_{2}{\bf S}\wedge{\bf S}_{y}+
a_{3}{\bf S}_{y}+(a_{3y}-\phi_{y}){\bf S}.
\eqno(20b)
$$
The compatibility condition of the equations (20) or the system (11)  gives 
the stationary M-XIIIA equation
$$
{\bf S}\wedge[a_{2}{\bf S}_{yy}+
(a_{1}-b_{2}){\bf S}_{xy}-b_{1}{\bf S}_{xx}]+
\phi_{y}{\bf S}_{x}+
\phi_{x}{\bf S}_{y}=0
\eqno(21a)
$$
$$
\phi_{xy}=\frac{(a_{1}+b_{2})}{2}{\bf S}\cdot ({\bf S}_{x}\wedge{\bf S}_{y}).
\eqno(21b)
$$
It follows from the M-XIIIA equation [28]
$$
{\bf S}_{t}={\bf S}\wedge[a_{2}{\bf S}_{yy}+
(a_{1}-b_{2}){\bf S}_{xy}-b_{1}{\bf S}_{xx}]+
\phi_{y}{\bf S}_{x}+
\phi_{x}{\bf S}_{y}
\eqno(22a)
$$
$$
\phi_{xy}=\frac{(a_{1}+b_{2})}{2}{\bf S}\cdot ({\bf S}_{x}\wedge{\bf S}_{y})
\eqno(22b)
$$
in the stationary limit.

\subsection{Stationary M-XIIIB equation}
Let 
$$
a_{5}=a_{3x}+\phi_{y}, \quad b_{5}=a_{3y}-\phi_{x} \eqno(23)
$$
so that instead of (8) we get the following MF [28]
$$
{\bf r}_{x}=a_{1}{\bf S}\wedge{\bf S}_{x}+a_{2}{\bf S}\wedge{\bf S}_{y}+
a_{3}{\bf S}_{x}+(a_{3x}+\phi_{y}){\bf S}
\eqno(24a)
$$
$$
{\bf r}_{y}=b_{1}{\bf S}\wedge{\bf S}_{x}+b_{2}{\bf S}\wedge{\bf S}_{y}+
a_{3}{\bf S}_{y}+(a_{3y}-\phi_{x}){\bf S}.
\eqno(24b)
$$
In this case  the system (11) takes the form
$$
{\bf S}\wedge[a_{2}{\bf S}_{yy}+
(a_{1}-b_{2}){\bf S}_{xy}-b_{1}{\bf S}_{xx}]+
\phi_{x}{\bf S}_{x}+
\phi_{y}{\bf S}_{y}=0
\eqno(25a)
$$
$$
\phi_{xx}+\phi_{yy}=(a_{1}+b_{2}){\bf S}\cdot ({\bf S}_{x}\wedge{\bf S}_{y})
\eqno(25b)
$$
that is the  stationary M-XIIIB
 equation. Itself the M-XIIIB equation can
be written as [28]
$$
{\bf S}_{t}={\bf S}\wedge[a_{2}{\bf S}_{yy}+
(a_{1}-b_{2}){\bf S}_{xy}-b_{1}{\bf S}_{xx}]+
\phi_{x}{\bf S}_{x}+
\phi_{y}{\bf S}_{y}
\eqno(26a)
$$
$$
\phi_{xx}+\phi_{yy}=(a_{1}+b_{2}){\bf S}\cdot ({\bf S}_{x}\wedge{\bf S}_{y}).
\eqno(26b)
$$

\subsection{Stationary Ishimori equation}
Similarly can be shown that the stationary Ishimori  equation
$$
{\bf S}\wedge({\bf S}_{xx}+\alpha^{2}{\bf S}_{yy})+
\phi_{x}{\bf S}_{y}+
\phi_{y}{\bf S}_{x}=0
\eqno(27a)
$$
$$
\alpha^{2}\phi_{yy}-\phi_{xx}=\alpha^{2}{\bf S}\cdot ({\bf S}_{x}\wedge{\bf S}_{y})
\eqno(27b)
$$
also can be obtained from the compatibility condition of the particular
case of the MF (8) [28].

\section{Surfaces induced by the generalized Myrzakulov formula}
Now we consider the following generalized MF [28]
$$
{\bf r}_{x}=a_{1}{\bf N}\wedge{\bf N}_{x}+a_{2}{\bf N}\wedge{\bf N}_{y}+
a_{3}{\bf N}_{x}+a_{4}{\bf N}_{y}+a_{5}{\bf N}+{\bf A}+C{\bf r}
+{\bf R}_{1}\wedge{\bf r}
\eqno(28a)
$$
$$
{\bf r}_{y}=b_{1}{\bf N}\wedge{\bf N}_{x}+b_{2}{\bf N}\wedge{\bf N}_{y}+
b_{3}{\bf N}_{x}+b_{4}{\bf N}_{y}+b_{5}{\bf N}+{\bf B}+D{\bf r}
+{\bf R}_{2}\wedge{\bf r}.
\eqno(28b)
$$
We cite a few particular examples of this generalized MF (for details see,
e.g., the ref. [28]).
\subsection{Example 1.}
Let us assume that the tangent vectors are given by
$$
{\bf r}_{x}=C {\bf r},\quad
{\bf r}_{y}=D {\bf r}\eqno(29)
$$
where $C,D$ are $3\times 3$ matrices. Hence we get
$$
C_{y}-D_{x}+[C,D]=0. \eqno(30)
$$
Let
$$
C=\left(
\begin{array}{ccc}
0&k&0\\
-k&0&\tau\\
0&-\tau&0
\end{array}
\right), \quad
D=\left(
\begin{array}{ccc}
0&\omega_{3}&-\omega_{2}\\
-\omega_{3}&0&\omega_{1}\\
\omega_{2}&\omega_{1}&0
\end{array}
\right).
\eqno(31)
$$
If we  introduce the complex function 
$
\psi=\frac{k}{2}e^{-i\partial^{-1}_{x} \tau}.$
Then from (30) we get [1]
$$
i\psi_{t}+\psi_{xx}+2|\psi|^{2}\psi=0 \eqno(32)
$$
which is the NLSE.
\subsection{Example 2.}
Now we consider the case when the tangent vectors are given by
$$
{\bf r}_{x}={\bf R}_{1}\wedge{\bf r},\quad
{\bf r}_{y}={\bf R}_{2}\wedge{\bf r}\eqno(33)
$$
where ${\bf R}_{i}$ are some vectors.  Hence we get
$$
{\bf R}_{1y}-{\bf R}_{2x}+2{\bf R}_{1}\wedge{\bf R}_{2}=0. \eqno(34)
$$

\section{Conclusion}

In this note, we have considered some 
 spin surfaces induced by SS in 
two dimensions. 
The corresponding SS are presented. These SS include,
the stationary Ishimori equation, Heisenberg ferromagnet, the
LLE and so on.  At last, in Appendix, we presented some generalized SS,
so-called Myrzakulov equations, which describe spin-phonon or magnetoelastic
systems.
\section{\it Appendix. Some magnetoelastic systems in 1+1 dimensions}
In this Appendix we present some spin-phonon or magnetoelastic systems
in 1+1 dimensions, so-called Myrzakulov equations, 
which were obtained in [28]. Some of these equations are integrable,
for instance, the Myrzakulov V (M-V), Myrzakulov XXXIV (M-XXXIV) and Myrzakulov LXIX
(LXIX) equations. 
\\
\vspace{0.3cm}
%\newpage 
\\
\hfill Table 1.  {\it Nonlinear magnetoelastic systems of 
the 0-type (The Landau-Lifshitz equations with potentials)}\\
\begin{tabular}{|l|c|}                  \hline
Name of the equation & Equation of motion\\ \hline
the  M-LVII equation &
$
2iS_t=[S,S_{xx}]+u[S,\sigma_3] $ \\ \hline
the M-LVI equation  &
$
2iS_t=[S,S_{xx}]+uS_3[S,\sigma_3] $ \\ \hline
the  M-LV equation  &
$
2iS_t=\{(\mu \vec S^2_x-u+m)[S,S_x]\}_x $ \\ \hline
the  M-LIV equation  &
$
2iS_t=n[S,S_{xxxx}]+2\{(\mu \vec S^2_x-u+m)[S,S_x]\}_x
 $ \\ \hline
the  M-LIII equation &
$
2iS_t=[S,S_{xx}]+2iuS_x $ \\ \hline
\end{tabular}
\\
\vspace{0.3cm}
\\
\hfill Table 2.  {\it Nonlinear magnetoelastic systems of 
the 1-type}\\
\begin{tabular}{|l|c|}                  \hline
Name of the equation & Equation of motion\\ \hline
the  M-LII equation &
$
2iS_t=[S,S_{xx}]+u[S,\sigma_3] 
$\\
&
$
\rho u_{tt}=\nu^2_0 u_{xx}+\lambda(S_3)_{xx}     
$ \\  \hline
the  M-LI equation &
$
2iS_t=[S,S_{xx}]+u[S,\sigma_3] 
$ \\ &
$
\hspace{2.5cm}\rho u_{tt}=\nu^2_0 u_{xx}+\alpha(u^2)_{xx}+\beta u_{xxxx}+
    \lambda(S_3)_{xx}                            
$ \\ \hline
the  M-L equation & 
$
2iS_t=[S,S_{xx}]+u[S,\sigma_3] 
$ \\ &
$
u_t+u_x+\lambda(S_3)_x=0 
$ \\ \hline
the  M-XLIX equation&
$
2iS_t=[S,S_{xx}]+u[S,\sigma_3] 
$ \\ &
$
u_t+u_x+\alpha(u^2)_x+\beta u_{xxx}+\lambda(S_3)_x=0
$ \\ \hline
\end{tabular}
\\
\vspace{0.3cm}
\\
\hfill Table 3. {\it Nonlinear magnetoelastic systems of 
the 2-type}\\
\begin{tabular}{|l|c|}                  \hline
Name of the equation & Equation of motion\\ \hline
the  M-XLVIII equation &
$
2iS_t=[S,S_{xx}]+uS_3[S,\sigma_3] 
$ \\ &
$
\rho u_{tt}=\nu^2_0 u_{xx}+\lambda(S^2_3)_{xx}     $ \\ \hline
the  M-XLVII equation &
$
2iS_t=[S,S_{xx}]+uS_3[S,\sigma_3] 
$ \\ & 
$
\hspace{2.5cm}\rho u_{tt}=\nu^2_0 u_{xx}+\alpha(u^2)_{xx}+\beta u_{xxxx}+
\lambda (S^2_3)_{xx} 
$ \\ \hline
the M-XLVI equation &
$
2iS_t=[S,S_{xx}]+uS_3[S,\sigma_3] 
$ \\ &
$
u_t+u_x+\lambda(S^2_3)_x=0 
$ \\ \hline
the  M-XLV equation &
$
2iS_t=[S,S_{xx}]+uS_3[S,\sigma_3] 
$ \\ &
$
u_t+u_x+\alpha(u^2)_x+\beta u_{xxx}+\lambda(S^2_3)_x=0
$\\ \hline
\end{tabular}
\\
\vspace{0.3cm}
%\newpage
\\
\hfill Table 4. {\it Nonlinear magnetoelastic systems of 
the 3-type}\\
\begin{tabular}{|l|c|}                  \hline
Name of the equation & Equation of motion\\ \hline
the  M-XLIV equation&
$
2iS_t=\{(\mu \vec S^2_x - u +m)[S,S_x]\}_x 
$ \\ &
$
\rho u _{tt}=\nu^2_0 u_{xx}+\lambda(\vec S^2_x)_{xx}   
$ \\ \hline
the  M-XLIII equation &
$ 
2iS_t=\{(\mu \vec S^2_x - u +m)[S,S_x]\}_x
$ \\ &
$
\hspace{2.5cm}\rho u _{tt}=\nu^2_0 u_{xx}+\alpha (u^2)_{xx}+\beta u_{xxxx}+ \lambda
(\vec S^2_x)_{xx} 
$ \\ \hline
the  M-XLII equation &
$
2iS_t=\{(\mu \vec S^2_x - u +m)[S,S_x]\}_x 
$ \\ &
$
u_t+u_x +\lambda (\vec S^2_x)_x = 0                     
$ \\ \hline
the  M-XLI equation &
$
2iS_t=\{(\mu \vec S^2_x - u +m)[S,S_x]\}_x 
$ \\ &
$
u_t+u_x +\alpha(u^2)_x+\beta u_{xxx}+\lambda (\vec S^2_x)_{x} = 0
$ \\ \hline
\end{tabular}
\\
\vspace{0.3cm}
\\
\hfill Table 5. {\it Nonlinear magnetoelastic systems of 
the 4-type}\\
\begin{tabular}{|l|c|} \hline
Name of the equation & Equation of motion\\ \hline
the  M-XL equation&
$
\hspace{0.2cm} 2iS_t=[S,S_{xxxx}]+2\{(\mu\vec S^2_x-u+m)[S,S_x]\}_{x} 
$ \\ &
$
\rho u_{tt}=\nu^2_0 u_{xx}+\lambda (\vec S^2_x)_{xx}  
$ \\ \hline
the  M-XXXIX  equation &
$
2iS_t=[S,S_{xxxx}]+2\{(\mu\vec S^2_x-u+m)[S,S_x]\}_{x}
$ \\ &
$
\rho u_{tt}=\nu^2_0 u_{xx}+\alpha(u^2)_{xx}+\beta u_{xxxx}+\lambda
(\vec S^2_x)_{xx} 
$ \\ \hline
the  M-XXXVIII equation &
$
2iS_t=[S,S_{xxxx}]+2\{(\mu\vec S^2_x-u+m)[S,S_x]\}_{x}
$\\ &
$
u_t + u_x + \lambda (\vec S^2_x)_x = 0
$ \\ \hline
the  M-XXXVII equation &
$
2iS_t=[S,S_{xxxx}]+2\{(\mu\vec S^2_x-u+m)[S,S_x]\}_{x} 
$ \\ &
$
u_t + u_x + \alpha(u^2)_x + \beta u_{xxx}+\lambda (\vec S^2_x)_x = 0
$ \\ \hline
\end{tabular}
\\
\vspace{0.3cm}
\\
\hfill Table 6. {\it Nonlinear magnetoelastic systems of 
the 5-type}\\
\begin{tabular}{|l|c|} \hline
Name of the equation & Equation of motion\\ \hline
the  M-XXXVI equation &
$
2iS_t=[S,S_{xx}]+2iuS_x   
$ \\ &
$
\rho u_{tt}=\nu^2_0 u_{xx}+\frac{\lambda}{4}(tr(S^2_x))_{xx} 
$ \\ \hline
the  M-XXXV equation &
$
2iS_t=[S,S_{xx}]+2iuS_x  
$ \\ &
$
\hspace{1.3cm}\rho u_{tt}=\nu^2_0 u_{xx}+\alpha(u^2)_{xx}+\beta u_{xxxx}+
\frac{\lambda}{4}(tr(S^2_x))_{xx} 
$ \\ \hline
the  M-XXXIV equation &
$
2iS_t=[S,S_{xx}]+2iuS_x
$\\ &
$
u_t + u_x + \frac{\lambda}{4}(tr(S^2_x))_{x} = 0 
$ \\ \hline
the  M-XXXIII equation &
$
2iS_t=[S,S_{xx}]+2iuS_x  
$ \\ &
$
u_t + u_x + \alpha(u^2)_x + \beta u_{xxx}+\frac{\lambda}{4}(tr(S^2_x))_{x} = 0
$\\ \hline
\end{tabular}
\\
\vspace{0.3cm}
%\newpage
\\
\hfill Table 7. {\it Nonlinear magnetoelastic systems of 
the 6-type}\\
\begin{tabular}{|l|c|} \hline
Name of the equation & Equation of motion\\ \hline
the  M-LXIX equation &
$
\hspace{2.9cm}{\bf S}_t=\frac{1}{\sqrt {{\bf S}^2_x}}(-\sqrt {{\bf S}^2_x-u^2}{\bf S}_x+u{\bf
S}\wedge{\bf S}_x)\hspace{1cm}   
$ \\ &
$
u_{x}=v\sqrt {{\bf S}^2_t-u^2} 
$ \\ &
$ v_t=-{\bf S}\cdot ({\bf S}_t\wedge{\bf S}_x)
$ \\ \hline
\end{tabular}
\\
\vspace{0.3cm}
%\newpage
\\
\hfill Table 8. {\it Nonlinear magnetoelastic systems of 
the 7-type}\\
\begin{tabular}{|l|c|} \hline
Name of the equation & Equation of motion\\ \hline
the  M-V equation &
$
\hspace{1cm}iS_{t} = \frac{1}{2}[S, S_{xx}] + \frac{3}{2}[S^{2}, (S^{2})_{xx}],\quad
S^{3}=S \in osp(2|1)\hspace{0.7cm}
$ \\ \hline
\end{tabular}

\end{document}